\documentclass[reqno,10pt]{amsart}

\usepackage{amsmath,amssymb,amsthm,amsfonts}
\usepackage{enumerate}
\usepackage{a4wide}
\newcommand{\invhom}[2]{\varphi_{{#2}}^{{#1}}}
\newcommand{\init}{\mathrm{in}}
\newcommand{\Tor}{\mathrm{Tor}}
\newcommand{\primes}{\mathcal{P}}

\theoremstyle{definition}
\newtheorem{definition}{Definition}[section]

\newtheorem{example}[definition]{Example}
\theoremstyle{remark}

\theoremstyle{plain}
\newtheorem{lemma}[definition]{Lemma}

\newtheorem{proposition}[definition]{Proposition}
\newtheorem{theorem}[definition]{Theorem}
\newtheorem{thm}[definition]{Theorem}
\newtheorem{corr}[definition]{Corollary}

\newtheorem{conjecture}[definition]{Conjecture}

\newcommand{\Nat}{{\mathbb{N}}}

\def\setsuchas#1#2{\left\{\,{#1}\,\vrule\,{#2}\,\right\}}
\newcommand{\set}[1]{{\{#1\}}}

\newcommand{\M}{\mathcal{M}}
\newcommand{\tdeg}[1]{\vert#1\rvert}
\newcommand{\supp}{\operatorname*{Supp}}

\newcommand{\divides}[2]{{#1 \left\lvert {#2} \right.}}

\begin{document}
\date{\today}
\title{Truncations of the ring of 
  number-theoretic functions} 
\author{Jan Snellman}

\date{\today}
\address{School of Informatics\\ 
  University of Wales\\
  Dean Street, Bangor \\
  Gwynedd LL57 1UT\\
  Wales\\
  UK
}

\email{jans@matematik.su.se}
\thanks{The author was supported by grants from Svenska Institutet and
  Kungliga Vetenskapsakademin, 
  and by grant n. 231801F from Centre International des Etudiants et
  Stagiaires.}

\begin{abstract}
  We study  the ring \(\Gamma\) of all
  functions \(\Nat^+ \to K\), endowed with the usual convolution
  product.  \(\Gamma\), which we call \emph{the ring of
  number-theoretic functions},  is an inverse limit of the
  ``truncations''  
  \[\Gamma_n =  \setsuchas{f \in \Gamma}{\forall m 
    > n: \, f(m)=0}.\]
  Each \(\Gamma_n\) is a zero-dimensional,
  finitely generated \(K\)-algebra, which may be expressed as the
  quotient of a finitely generated polynomial ring with a
  \emph{stable} (after reversing the order of the variables) monomial
  ideal. Using the description of the 
  free minimal resolution 
  of stable ideals given by Eliahou-Kervaire, and some additional
  arguments by Aramova-Herzog and Peeva, we give the
  Poincar\'e-Betti  series for \(\Gamma_n\).
\end{abstract}

\maketitle

\begin{section}{Introduction}
  Cashwell and Everett \cite{NumThe} studied ``the ring of
  number-theoretic functions'' 
  \begin{equation}
    \label{eq:gammadef}
    \Gamma = \setsuchas{f}{\Nat^+ \to K}
  \end{equation}
  where \(\Nat^+\) is the set of positive natural numbers (we denote
  by \(\Nat\) the set of all natural numbers) and \(K\) is a field
  containing the rational numbers. \(\Gamma\) is endowed with
  component-wise addition and multiplication with scalars, and with
  the convolution (or Cauchy) product
  \begin{equation}\label{eq:convdef}
    fg(n) = \sum_{\substack{(a,b) \in {(\Nat^+)} \times {(\Nat^+)}
    \\ab=n}} 
    f(a)g(b) 
  \end{equation}
  With these operations, \(\Gamma\) becomes a commutative
  \(K\)-algebra. It is immediate that it is a local domain; less
  obvious is the fact that it is a unique factorisation domain. 
  Cashwell and Everett proved this in \cite{NumThe} using the
  isomorphism
  \begin{equation}\label{eqn:gpiso}
    \begin{split}
      \Phi: \Gamma & \to  K[[X]] \\
      f & \mapsto \sum f(n) x_1^{\alpha_1}x_2^{\alpha_2}\cdots
    \end{split}
  \end{equation}
  where \(X=\set{x_1,x_2,x_3,\dots}\), \(K[[X]]\) is the  ``large''
  power series ring of all functions from the free abelian monoid
  \(\M = [X] \) (the free abelian monoid generated by \(X\)) to \(K\),
  and where  
  the summation extends over all
  \(n=p_1^{\alpha_1} p_2^{\alpha_2}\cdots \in \Nat^+\). Here, and
  henceforth, we denote by \(p_i\) the \(i\)'th prime
  number, with \(p_1=2\), and by \(\primes\) the set of all prime
  numbers. That \eqref{eqn:gpiso} is an isomorphism is immediate
  from the following isomorphism of commutative monoids, implied by 
  the fundamental theorem of arithmetics:
  \begin{equation}\label{eqn:fundiso}
    (\Nat^+,\cdot) \simeq \coprod_{p \in \primes}{(\Nat,+)}
  \end{equation}

  The following
  number-theoretic functions are of particular interest (whenever
  possible, we use the same 
  notation as in \cite{NumThe}): 
  \begin{enumerate}
  \item The multiplicative unit \(\epsilon\) given by
    \(\epsilon(1)=1\), \(\epsilon(n) = 0\) for \(n > 1\),
  \item
    \(\lambda: \Nat^+ \to \Nat\) given by
    \(\lambda(1)=0\), \(\lambda(q_1\cdots q_l)=l\) if \(q_1,\dots,q_l\)
    are any (not necessarily distinct) prime numbers.
  \item     \(\tilde{\lambda}: \Nat^+ \to \Nat\) given \(\tilde{\lambda}(1)=0\),
    \(\tilde{\lambda}(p_1^{a_1} \dots p_r^{a_r}) = \sum a_r p_r\). 
  \item The M\"obius function \(\mu(1) = 1\), \(\mu(n) = (-1)^v\) if
    \(n\) is the product of \(v\) distinct prime factors, and \(0\)
    otherwise, 
  \item   For any \(i \in \Nat^+\),
  \(\chi_i(p_i) = 1\), and \(\chi_i(m)=0\) for \(m \neq p_i\). 
  Note that under the isomorphism \eqref{eqn:gpiso},
  \(\Phi(\chi_i)=x_i\). 
  \end{enumerate}

  The topic of this article is the study of the ``truncations''
  \(\Gamma_n\), where for each \(n \in \Nat^+\),
  \begin{equation}\label{eq:defgamman}
    \Gamma_n = \setsuchas{f \in \Gamma}{m > n \implies f(m) = 0}
  \end{equation}
  With the modified multiplication given by 
  \begin{equation}
    fg(n) = \sum_{\substack{(a,b) \in {\set{1,\dots,n}} \times {\set{1,\dots,n}}
    \\ab=n}} 
    f(a)g(b) 
  \end{equation}
  \(\Gamma_n\) becomes a \(K\)-algebra, isomorphic to \(\Gamma/J_n\), where
  \(J_n\) is the ideal 
  \[J_n=\setsuchas{f \in \Gamma}{\forall m \le n: \,
    f(m)=0}.\]
  If we define
  \begin{align}\label{eq:gammamap}
    \pi_n: \Gamma &\to \Gamma_n \\
    \pi_n(f)(m)  & = 
    \begin{cases}
      f(m) & m \le n \\
      0 & m > n
    \end{cases}
  \end{align}
  then \(\pi_n\) is a \(K\)-algebra epimorphism, and \(J_n\) is the kernel of
  \(\pi_n\). 
  We note furthermore that \(J_n\) is
  generated by \emph{monomials} in the elements \(\chi_i\).
  
  To describe the main idea of this paper, we need a few additional
  definitions. First,  for any \(n \in \Nat^+\) we denote by \(r(n)
  \in \Nat\) the  largest integer such that \(p_{r(n)} \le n\). In
  other words, \(r(n)\)  is the number of prime numbers \(\le n\)
  (this number is often denoted \(\pi(n)\)).
  Secondly, for a
  monomial \(m = x_1^{\alpha_1}\cdots x_w^{\alpha_w}\), we define the
  \emph{support} \(\supp(m)\) as the set of positive 
  integers \(i\) such that \(\alpha_i > 0\). We define 
  \(\max(m)\) and \(\min(m)\) as the maximal and minimal elements in the
  support of \(m\).
  \begin{definition}\label{def:stab}
    A monomial ideal \(I \subset K[x_1,\dots,x_r]\) is said to be
    \emph{strongly 
    stable} if whenever \(m\) is a monomial such that \(x_j m \in I\),
    then \(x_i m \in I\) for all \(i \le j\). If this condition holds 
    at least for all \(i \le j=\max(m)\) then \(I\) is said to be
    \emph{stable}. 
  \end{definition}
  
  We can now state our main theorem:
  \begin{thm}\label{thm:main}
    Let \(n \in \Nat^+\) and \(r=r(n)\). Then the following holds:
    \begin{enumerate}[(I)]
    \item \label{enum:monquotient}    
      \(\Gamma_n \simeq \frac{K[x_1,\dots,x_r]}{I_n}\), where 
      \(I_n\) is  a strongly stable monomial ideal, with respect to the
      reverse order of the variables. 
    \item \label{enum:hilbert}    
      \(\Gamma_n\) is artinian, with
      \(\dim_K(\Gamma_n)=n\). Furthermore, if it is given the natural
      grading with \(\tdeg{\chi_i}=1\), then its Hilbert series is
     \(\sum_{i} d_{i} t^i\)
    where \(d_i\) is the number of \(w \le n\)
    with \(\lambda(w) = i\).
    \item  \label{enum:oneone}    
      There is a 1-1
    bijection between the minimal monomial generators of \(I_n\) of
    minimal support \(v\), and
    the solutions in non-negative integers to the equation
    \begin{equation}
      \label{eqn:fcnvd}
      \begin{split}
      \log n - \log p_v & < \sum_{i=v}^r b_i \log p_i \le \log n 
      \end{split}
    \end{equation}
  \item \label{enum:PB}    
    If we denote by \(C_{n,v}\) the number of such solutions, then
    the Poincar\'e-Betti series of the free minimal resolution of
    \(K\) as a cyclic module over \(\Gamma_n\) is the following
    rational function:
    \begin{equation}
      \label{eq:pbser}
      P(\Tor_{*}^{\Gamma_n}(K,K),t) =    
      \frac{(1+t)^r}{1- t^2 \left(\sum_{i=1}^r
          (1+t)^{(i-1)}C_{n,r-i+1}\right)} 
    \end{equation}
    \end{enumerate}
  \end{thm}
  
  We will show this result, and also give the graded Poincar\'e-Betti
  series. 
  For this, we define the number \(C_{n,v,d}\) which counts the number
  of minimal generators of \(I_n\) of minimal support \(v\) and total
  degree \(d\). We determine some elementary properties of the numbers
  \(C_{n,v,d}\) and \(C_{n,v}\).

\end{section}

\begin{section}{The ring of number-theoretic functions  and its
  truncations}   
\begin{subsection}{Norms, degrees, and multiplicativity}
For a monomial \(\M \ni m = x_1^{a_1}\dots x_n^{a_n}\) we define the
\emph{weight} of \(m\) as 
\(w(m) = p_1^{a_1}\dots p_n^{a_n}\) (we put \(w(1)=1\)). Hence \(w\)
gives a bijection 
between \(\M\) and \(\Nat^+\).
Furthermore, we can define a term order on \(\M\) by \(m > m'\) iff
\(w(m) > w(m')\). If we define the \emph{initial monomial}
\(\init(f)\) 
of \(f \in K[[X]]\) as the monomial in \(\supp(f)\) minimal with
respect to \(>\), then \(\init(f)\) is easily seen to correspond to the
\emph{norm} \(N(\alpha)\) of a number-theoretic function \(\alpha\),
defined as the 
smallest \(n\) such that \(\alpha(n) \neq 0\). Here, we must use \(w\)
and \(\Phi\) to identify \(\M\) and \(\Nat^+\) and \(K[[X]]\) and
\(\Gamma\). As observed in \cite{NumThe}, the norm is multiplicative:
\(N(\alpha \beta) = N(\alpha) N(\beta)\).

Cashwell and Everett also define the \emph{degree} \(D(\alpha)\) to
mean the smallest \(d\) such that there exists an \(n\) with
\(\lambda(n)=d\) and \(\alpha(n) \neq 0\). This corresponds the
smallest \emph{total degree} of a monomial in
\(\supp(f)\). Furthermore, 
the norm \(M(\alpha)\), defined as the smallest integer \(n\) with
\(\lambda(n)=D(\alpha)\), \(\alpha(n) \neq 0\), corresponds to the
initial monomial of \(f\) under the term order obtained by refining
the total degree partial order with the term order \(>\).

A \emph{multiplicative function} is an element \(\alpha \in \Gamma\)
such 
that \(\alpha(1)=1\) and \(\alpha(ab) = \alpha(a)\alpha(b)\) whenever
\(a\) and \(b\) are 
relatively prime. Cashwell and Everett observes that a multiplicative
function is necessarily a unit in \(\Gamma\). One can further observe that
if \(\alpha\) is multiplicative, then \(f=\Phi(\alpha)\) can be
written
\begin{displaymath}
  f(x_1,x_2,x_3,\dots) = f_1(x_1) f_2(x_2) f_3(x_3) \cdots
\end{displaymath}
where each \(f_i(x_i) \in K[[x_i]]\) is invertible. In particular, the
constant function \(\Gamma \ni \nu_0\) with \(\nu_0(n)=1\) for all
\(n\), corresponds to
\begin{displaymath}
  \sum_{m \in \M} m = \frac{1}{1-x_1}\frac{1}{1-x_2}\frac{1}{1-x_3}
  \cdots 
\end{displaymath}
Since the M\"obius function is defined to be the inverse of this
function, we get that it corresponds to
\begin{displaymath}
  (1-x_1)(1-x_2)(1-x_3)\cdots = 1 - (\sum_{i=1}^\infty x_i) + 
  (\sum_{i < j} x_i x_j) - (\sum_{i < j < k} x_i x_j x_k) + \cdots
\end{displaymath}
\end{subsection}

\begin{subsection}{Truncations of the ring of number-theoretic
  functions} 
Let \(n,n' \in \Nat^+\), \(n' > n\). Then there is a \(K\)-algebra epimorphism 
\begin{align*}
\invhom{n'}{n}: \Gamma_{n'} &\to \Gamma_n \\
\invhom{n'}{n}(f)(m)  & = 
\begin{cases}
  f(m) & m \le n \\
  0 & m > n
\end{cases}
\end{align*}

Hence, the \(\Gamma_n\)'s form an inverse system.
\begin{lemma}\label{lemma:invlim}
  \(\varprojlim \Gamma_n \simeq \Gamma\).
\end{lemma}
\begin{proof}
Given any \(f \in \Gamma\), the sequence \((\pi_1(f), \pi_2(f),
\pi_3(f),\dots)\) is coherent. Conversely, given any coherent sequence
\((g_1,g_2,g_3,\dots)\), we can define \(g: \Nat \to K\) by \(g(m) =
g_i(m)\) where \(i \ge m\).
\end{proof}
As a side remark, we note that 
\begin{lemma}
  The decreasing filtration 
  \begin{equation}
    \label{eq:Jnfilt}
    J_1 \supsetneq J_2 \supsetneq J_3 \supsetneq \cdots
  \end{equation}
  is separated, that is, \(\cap_{n} J_n = (0)\).
\end{lemma}

\begin{definition}
We define 
\begin{equation}\label{eqn:In}
  I_n = K[[X]] \setsuchas{m \in \M}{w(m) > n},
\end{equation}
that is, as the monomial ideal in \(K[[X]]\) generated by all
monomials of weight strictly higher than \(n\). 
We put  \(A_n = \frac{K[[X]]}{I_n}\).  
\end{definition}

\begin{proposition}
A \(K\)-basis of \(A_n\) is given by all monomials of weight \(\le n\).
Hence \(A_n\) is an artinian algebra, with \(\dim_K(A_n) = n\).
Putting \(r=r(n)\), we have that 
\begin{equation}\label{eq:An}
  A_n = \frac{K[[X]]}{I_n} \simeq \frac{K[x_1,\dots,x_r]}{I_n \cap
    K[x_1,\dots,x_r]} 
\end{equation}  
\end{proposition}
\begin{proof}
  As a vector space, \(K[[X]] \simeq U \oplus I_n\), where \(U\)
  consists of all functions supported on monomials of weight 
  \(\le  n\). It follows that \(A_n \simeq U\) as \(K\) vector
  spaces. Of course, there are exactly \(n\) monomials of weight 
  \(\le n\). Finally, if \(s > r\) then \(w(x_s) = p_s > n\), hence
  \(x_s \in I_n\).
\end{proof}

We will abuse notations and identify \(I_n\) and its contraction 
\(I_n \cap K[x_1,\dots,x_r]\).

\begin{lemma}\label{lemma:trunciso}
  \(\Gamma_n \simeq A_n\).
\end{lemma}
\begin{proof}
  Since \(A_n\) has a \(K\)-basis is given by all monomials of
  weight \(\le n\), the two \(K\)-algebras are isomorphic as
  \(K\)-vector spaces. The multiplication in \(A_n\) is induced from
  the multiplication in \(K[[X]]\), with the extra condition that
  monomials of weight \(>n\) are truncated. This is the same
  multiplication as in \(\Gamma_n\).
\end{proof}

  \begin{proposition}\label{prop:antistable}
    \(I_n\) is a strongly stable ideal, with respect to the reverse order of
    the variables.
  \end{proposition}
  \begin{proof}
    We must show that if \(m \in I_n\), and \(\divides{x_i}{m}\), then 
    \(m x_j/x_i \in I\) for \(i \le j \le r\).
    We have that \(w(m x_j/x_i) = w(m) p_j/p_i > w(m) > n\).
  \end{proof}
  Part~\ref{enum:monquotient} of the main theorem is now proved.
  

  We give \(K[x_1,\dots,x_r]\) an \(\Nat^2\)-grading by giving the
  variable \(x_i\) bi-degree \((1,p_i)\). Since each \(I_n\) is
  bihomogeneous, this grading is inherited by \(A_n\).
  \begin{theorem}\label{thm:hilbert}
    The bi-graded Hilbert series of \(A_n\) is given by
    \begin{displaymath}
      A_n(t,u) = \sum_{i,j} c_{ij} t^i u^j, 
    \end{displaymath}
    where \(c_{ij}\) is the number of \(p_1^{a_1} \dots p_r^{a_r} 
    \le n\) with \(\sum a_r = i\) and \(\sum a_r p_r = j\).
    Furthermore, 
    \begin{align*}
      A_n(t,1) &= \sum_{i} d_{i} t^i \\
      A_n(1,u) &= \sum_{j} e_{j} u^j 
    \end{align*}
    where \(d_i\) is the number of \(w \le n\)
    with \(\lambda(w) = i\), and  
    \(e_i\) is the number of \(w \le n\)
    with \(\tilde{\lambda}(w) = i\). In particular, the 
    \(t^1\)-coefficient of \(A_n(t,1)\) is the number of prime numbers
    \(\le n\).
  \end{theorem}
  \begin{proof}
    The monomial \(x_1^{a_1}\cdots x_n^{a_n}\) has bi-degree
    \((\sum_{i=1}^n a_i, \sum a_i p_i)\).
  \end{proof}
  This establishes part~\ref{enum:hilbert} of the main theorem.   
\end{subsection}
\end{section}

\begin{section}{Minimal generators for \(I_n\)}
  Let \(n \in \Nat^+\), and let \(r=r(n)\).
  We have that
  \begin{equation}
    \label{eqn:trivweight}
    x_1^{a_1} \dots x_r^{a_r} = m \in I_n \quad \iff \quad w(m) > n
    \quad \iff \quad \prod_{i=1}^r p_i^{a_i} > n. 
  \end{equation}
  We denote by \(G(I_n)\) the set of minimal monomial generators of
  \(I_n\). 
  For \(m = x_1^{a_1} \dots x_r^{a_r}\) to be an element of \(G(I_n)\)
  it is necessary and 
  sufficient that \(m \in I_n\) and that for \(1 \le v \le r\), 
  \(\divides{x_v}{m} \implies m/x_v \not \in I_n\). In other words,
  \begin{equation}
    \label{eqn:mingen}
    1 \le j \le n, \, a_j > 0 \quad \implies \quad n < \prod_{i=1}^r
    p_i^{a_i} \le p_j n. 
  \end{equation}

  \begin{definition}
    For \(n,v,d\) positive integers, we define:
    \begin{align}
      C_n &= \#G(I_n) \\ 
    C_{n,v} &= \#\setsuchas{m \in G(I_n)}{\min(m)=v} \\
    C_{n,v,d} &= \#\setsuchas{m \in G(I_n)}{\min(m)=v, \, \tdeg{m}=d}
    \end{align}
  \end{definition}

  \begin{figure}[p]
    \begin{center}
      \leavevmode
      \caption{The numbers \(C_n\) and \(C_{n,i}\).}
      \label{fig:Cni}
      \tiny{
      \begin{displaymath}
\begin {array}{|c|c|c|c|c|c|c|c|c|c|c|c|} 
\hline
n & \Sigma & i=1 & i=2 & 3 & 4 & 5 & 6 & 7 & 8 &9 & 10\\ \hline
2&1&1
\\3&3&2&1
\\4&3&2&1
\\5&6&3&2&1
\\6&6&3&2&1
\\7&10&4&3&2&1
\\8&10&4&3&2&1
\\9&11&5&3&2&1
\\10&11&5&3&2&1
\\11&16&6&4&3&2&1
\\12&16&6&4&3&2&1
\\13&22&7&5&4&3&2&1
\\14&22&7&5&4&3&2&1
\\15&23&8&5&4&3&2&1
\\16&23&8&5&4&3&2&1
\\17&30&9&6&5&4&3&2&1
\\18&30&9&6&5&4&3&2&1
\\19&38&10&7&6&5&4&3&2&1
\\20&38&10&7&6&5&4&3&2&1 \\
21 & 39 & 11 & 7 & 6 & 5  & 4 & 3 & 2 & 1\\
22 & 39 & 11 & 7 & 6 & 5  & 4 & 3 & 2 & 1\\
23 & 48 & 12 & 8 & 7 & 6 & 5 & 4 & 3 & 2 & 1 \\
24 & 48 & 12 & 8 & 7 & 6 & 5 & 4 & 3 & 2 & 1 \\
25 & 50 & 13 & 9 & 7 & 6 & 5 & 4 & 3 & 2 & 1 \\
26 & 50 & 13 & 9 & 7 & 6 & 5 & 4 & 3 & 2 & 1 \\
27 & 51 & 14& 9 & 7 & 6 & 5 & 4 & 3 & 2 & 1 \\
28 & 51 & 14& 9 & 7 & 6 & 5 & 4 & 3 & 2 & 1 \\
29 & 61 & 15 & 10 & 8 & 7 & 6 & 5 & 4 & 3 & 2 & 1 \\
30 & 61 & 15 & 10 & 8 & 7 & 6 & 5 & 4 & 3 & 2 & 1 \\ \hline
\end {array}
      \end{displaymath}

}
    \end{center}
  \end{figure}

  \begin{figure}[p]
    \begin{center}
      \leavevmode
      \caption{The numbers \(C_{n,i,g}\).}
      \label{fig:Cnig}
      \tiny{
      \begin{displaymath}
        \begin {array}{|c|c|c|c|c|c|c|c|c|c|c|c|} 
          \hline
          n & i=1 & i=2 & 3 & 4 & 5 & 6 & 7 &8 & 9\\ \hline
2&1\\
3&2&1\\
4&u+1&1\\
5&u+2&2&1\\
6&2\,u+1&2&1\\
7&2\,u+2&3&2&1\\
8&{u}^{2}+u+2&3&2&1\\
9&{u}^{2}+2\,u+2&u+2&2&1\\
10&{u}^{2}+3\,u+1&u+2&2&1\\
11&{u}^{2}+3\,u+2&u+3&3&2&1\\
12&2\,{u}^{2}+2\,u+2&u+3&3&2&1
\\13&2\,{u}^{2}+2\,u+3&u+4&4&3&2&1
\\14&2\,{u}^{2}+3\,u+2&u+4&4&3&2&1
\\15&2\,{u}^{2}+4\,u+2&2\,u+3&4&3&
2&1\\16&{u}^{3}+{u}^{2}+4\,u+2&2
\,u+3&4&3&2&1\\17&{u}^{3}+{u}^{2
}+4\,u+3&2\,u+4&5&4&3&2&1\\18&{u}^
{3}+2\,{u}^{2}+3\,u+3&2\,u+4&5&4&3&2&1\\
19&{u}^{3}+2\,{u}^{2}+3\,u+4&2\,u+5&6&
5&4&3&2&1\\
20&{u}^{3}+3\,{u}^{2}+2\,
u+4&2\,u+5&6&5&4&3&2&1\\
21&{u}^{3}+3\,{u}^{2}+3\,u+4&3\,u+4&6&5&4&3&2&1\\
22&{u}^{3}+3\,{u}^{2}+4\,u+3&3\,u+4&6&5&4&3&2&1\\23&{u}^{3}+3\,{u}^{2}+4\,
u+4&3\,u+5&7&6&5&4&3&2&1\\24&2\,{u}^{3
}+2\,{u}^{2}+4\,u+4&3\,u+5&7&6&5&4&3&2&1\\
25&2\,{u}^{3}+2\,{u}^{2}+5\,u+4&4\,u+5
&u+6&6&5&4&3&2&1 \\
26&2\,{u}^{3}+2\,{u}^{2}+6\,u+3&4\,u+5
&u+6&6&5&4&3&2&1\\
27&2\,{u}^{3}+3\,{u}^{2}+6\,u+3&{u}^{2}+3\,u+5&u+6&6&5&4&3&2&1
\\
28&2\,{u}^{3}+4\,{u}^{2}+5\,
u+3&{u}^{2}+3\,u+5&u+6&6&5&4&3&2&1
\\
29&2\,{u}^{3}+4\,{u}^{2}+5\,u+4&{u}^{2
}+3\,u+6&u+7&7&6&5&4&3&2&1\\
30&2\,{u}^{3}+5\,{u}^{2}+4\,u+4&{u}^{2}+3\,u+6&u+7&7&6&5&4
&3&2&1
\\ \hline
\end {array}
      \end{displaymath}

      }
    \end{center}
  \end{figure}

  \begin{theorem}\label{thm:cnv}
    \(C_{n,v}\) is the number of solutions \((b_1,\dots,b_r) \in
    \Nat^r\) to the equation
    \begin{equation}
      \label{eqn:cnv}
      \log n - \log p_v < \sum_{i=v}^r b_i \log p_i \le \log n.
    \end{equation}
    Equivalently, \(C_{n,v}\) is the number of integers \(x\) such that
    \(n/p_v < x \le n\) and such that no prime factors of \(x\) are smaller
    than \(p_v\). 

    Similarly,
    \(C_{n,v,d}\) is the number of solutions \((b_1,\dots,b_r) \in
    \Nat^r\) to the system of equations
    \begin{equation}
      \label{eqn:cnvd}
      \begin{split}
      \log n - \log p_v & < \sum_{i=v}^r b_i \log p_i \le \log n \\
      \sum_{i=1}^r b_i &= d-1.
      \end{split}
    \end{equation}
    or equivalently, \(C_{n,v,d}\) is the number of integers \(x\) such that
    \(n/p_v < x \le n\) and such that  no prime factors of \(x\) are smaller
    than \(p_v\), and with the additional constraint that
    \(\lambda(x)=d\). 
  \end{theorem}
  \begin{proof}
    We have that \(a_v >0\), \(a_w = 0\) for \(w < v\). Hence
    equation \eqref{eqn:mingen} implies that
    \begin{displaymath}
      \quad n < \prod_{j=v}^r p_i^{a_i} \le p_v n.
    \end{displaymath}
    Putting \(b_v = a_v - 1\), \(b_j = a_j\) for \(j > v\) we can
    write this as
    \begin{displaymath}
      \quad n < p_v \prod_{j=v}^r p_i^{b_i} \le p_v n \quad \iff \quad
      n/p_v <  \prod_{j=v}^r p_i^{b_i} \le n 
    \end{displaymath}
    from which \eqref{eqn:cnv} follows by taking logarithms. This
    implies \eqref{eqn:cnvd} as well.
  \end{proof}
  We have now proved  part~\ref{enum:oneone}     of the main theorem.
  \begin{example}
  The first few \(I_n\)'s are as follows: \(I_2 = (x_1^2)\),
  \(I_3 = (x_1^2,x_2^2,x_1x_2)\), \(I_4 =  (x_1^3,x_2^2,x_1x_2)\),
  \(I_5 = (x_1^3,x_2^2,x_1x_2,x_3^2,x_1x_3, x_2x_3)\).
  \end{example}
  We tabulate \(C_{n,i}\) and \(C_{n,i,d}\), the latter in form of the
  polynomial \(u^{-2}\sum_j C_{n,i,j} u^j\) in the tables \ref{fig:Cni} and
  \ref{fig:Cnig}.

  \begin{thm} \label{thm:bc}
    \begin{enumerate}[(1)]
    \item \label{enum:zero}
      \(C_{n,v}=0\) for \(v > r(n)\)
    \item \label{enum:geind}
      \(\forall n \in \Nat: \, \, \forall v \le r(n): \,
      C_{n,1+r(n)-v} \ge v\), 
    \item \label{enum:gesum}
      \(\forall n \in \Nat: \, C_n \ge \binom{r(n)+1}{2}\),
    \item \label{enum:lintail}
      \(\forall v \in \Nat: \, \exists N: \, \forall n \ge N:
      C_{n,1+r(n)-v} = v\).
    \item \label{enum:twostep}
      If \(n\) is even, then \(C_{n,v}=C_{n-1,v}\) for all \(v\), 
    \item \label{enum:ceil}
      \(C_{n,1} = \lceil n/2 \rceil\).
    \end{enumerate}
  \end{thm}
  \begin{proof}
    \eqref{enum:zero} Obvious. 
    
    \eqref{enum:geind}  and \eqref{enum:gesum}
    It suffices to show that for any subset \(S \subset
    \set{1,\dots,r}\) 
    of cardinality 1 or 2, there is an \(m \in G(I_n)\) with
    \(\supp(m)=S\). If \(S=\set{i}\) then there is an unique positive
    integer \(a\) such that \(p_i^{b-1} \le n < p_i^b\), and
    \(m=x_i^b\) is the desired generator. If \(S=\set{i,j}\) with \(i
    < j\) then we claim that there is a positive integer \(a\) such that
    \(x_i^a x_j \in G(I_n)\). Namely, choose \(b\) such that 
    \(p_i^{b-1} \leq n < p_i^b\), then since \(p_i < p_j\) one has \(n <
    p_i^{b-1} p_j\). Hence \(x_i^{b-1}x_j \in I_n\), so it is a multiple of
    some minimal generator. By the definition of \(b\), this minimal generator 
    must be of the form \(x_i^ax_j\) for some \(a\), which establishes the
    claim. 

    \eqref{enum:ceil} We must show that the number of
    solutions in \(\Nat^r\) to 
    \begin{displaymath}
      \frac{n}{2} < \prod_{i=1}^r {p_i}^{b_i} \le n
    \end{displaymath}
    is precisely \(\lceil \frac{n}{2} \rceil\). 
    Obviously, any integer \(\in (\frac{n}{2}, n]\) fits the bill;
    there are \(\lceil \frac{n}{2} \rceil\) of those.
    
    \eqref{enum:twostep} The case \(v=1\) follows from
    \eqref{enum:ceil}. Hence, it suffices to show that if \(v>1\), 
    \(x \in (\frac{n}{p_v}, n] \cap \Nat\), and if \(x\) has no prime
    factor \(< p_v\), then \(x \in (\frac{n-1}{p_v}, n-1] \cap
    \Nat\). The only way this can fail to happen is if \(x=n\), but
    then \(x\) is even, and has the prime factor \(2=p_1 < p_v\), a
    contradiction. 
    
    \eqref{enum:lintail}
    For large enough \(n\), the only integers \(x \le n\) with all
    prime factors \(\ge 1 + r(n) -v\) are
    \(p_{1+r(n)-v},\dots,p_{r(n)}\). There is \(v\) of these, and they
    are all \(> \frac{n}{p_v}\).
  \end{proof}

 \begin{thm} \label{thm:bcg}
    \begin{enumerate}
    \item \(C_{n,v,d}=0\) for \(v > r(n)\), and for \(d < 2\),
    \item \(\forall v \in \Nat: \, \exists N: \, \forall n \ge N:
      C_{n,1+r(n)-v,2} = v\), \(C_{n,1+r(n)-v,d} = 0\) for \(d \neq
      2\), 
    \item \(\binom{r(n)}{2} = \#\setsuchas{m \in \Nat^+}{m \le n, \,
        \lambda(m)=2}\). 
    \end{enumerate}
  \end{thm}
  \begin{proof}
    The first and the last assertions are obvious. The second one
    follows from the 
    proof of \eqref{enum:lintail} in the previous lemma. 
  \end{proof}

\end{section}

\begin{section}{Poincar{\'e} series}
  In     \cite{Eliahou:MinRes}, a minimal free multi-graded resolution
  of a \(I\) over \(S\) is given, where \(S=K[x_1,\dots,x_r]\) is a
  polynomial ring, and \(I \subset (x_1,\dots,x_r)^2\) is a stable
  ideal. As a consequence, the 
  following formula for the Poincar{\'e}-Betti series  is derived: 
  \begin{equation}\label{eqn:EliKer}
    P(\Tor_{*}^{S}(I,K),t) = \sum_{a \in G(I)} (1+t)^{\max(a) - 1}    
  \end{equation}
  where \(G(I)\) is the minimal generating set of \(I\).
  Since the resolution is multi-graded, \eqref{eqn:EliKer} can be
  modified to yield a formula for the graded Poincar{\'e}-Betti series
  (we here consider \(S\) as \(\Nat\)-graded, with each variable
  given weight 1):
  \begin{equation}\label{eqn:EliKerG}
    P(\Tor_{*,*}^{S}(I,K),t,u) = \sum_{a \in G(I)}
    u^{\tdeg{a}}(1+t)^{\max(a) - 1}     
  \end{equation}
  We will use the following variant of this result:
  \begin{theorem}[Eliahou-Kervaire]
    Let \(I \subset (x_1,\dots,x_r)^2 \subset K[x_1,\dots,x_r] = S\)
    be a stable monomial ideal. 
    Put 
      \begin{align}
        b_{i,d} &= \# \setsuchas{m \in G(I)}{\max(m)=i, \, \tdeg{m}=d}
        \\ 
        b_i &= \# \setsuchas{m \in G(I)}{\max(m)=i} 
      \end{align}
    Then 
    \begin{align}
        P(\Tor_{*}^{S}(I,K),t) &= \sum_{i=1}^r b_i (1+t)^{(i-1)} \\       
        P(\Tor_{*,*}^{S}(I,K),t,u) &= \sum_{i=1}^r \left((1+tu)^{(i-1)}
         \sum_j  b_{i,j} u^{j}  \right).       
  \end{align}
  For the Betti-numbers we have that
  \begin{equation}
    \beta_q = \dim_K \left(\Tor_q^S(I,K)\right) = \sum_{i=1}^r b_i
    \binom{i-1}{q}. 
  \end{equation}
\end{theorem}


From Proposition~\ref{prop:antistable} we have that the ideals \(I_n\)
are stable after reversing the order of
the variables. Hence, replacing \(\max\) by \(\min\), 
and hence \(b_i\) with \(C_{n,1+r-i}\), we get: 

\begin{corr}
  Let \(n \in \Nat^+\), \(r=r(n)\), \(S=K[x_1,\dots,x_r]\). Then
    \begin{align}
      P(\Tor_{*}^{S}(I_n,K),t) &= \sum_{i=1}^r C_{n,1+r-i} (1+t)^{(i-1)} \\
      P(\Tor_{*,*}^{S}(I_n,K),t,u) &= \sum_{i=1}^r (1+tu)^{(i-1)} \sum_j
      C_{n,1+r-i,j} u^{j}.       
    \end{align}
  For the Betti-numbers we have that
  \begin{equation}
    \beta_q = \sum_{i=1}^r C_{n,1+r-i} \binom{i-1}{q}.
  \end{equation}
\end{corr}

In \cite{Peeva:Stable, KoszulEliKer} it is shown that if
\(S=K[x_1,\dots,x_r]\) and \(I\) is a stable monomial ideal in \(S\),
then \(S/I\) is a Golod ring. Hence, from a
result of Golod \cite{Golod:Res} (see also \cite{GullLev}), it follows
that 
\begin{equation}
  \label{eqn:golod}
  P(\Tor_{*}^{S/I}(K,K),t) = \frac{(1+t)^r}{1- t^2 P(\Tor_{*}^{S}(I,K),t)}
\end{equation}
Regarding \(S\) as an \(\Nat\)-graded ring, one can show that in fact
\begin{equation}
  \label{eqn:golodGrad}
  P(\Tor_{*}^{S/I}(K,K),t,u) = \frac{(1+ut)^r}{1- t^2 P(\Tor_{*}^{S}(I,K),t,u)}
\end{equation}

The following theorem is an immediate consequence:
\begin{theorem}[Herzog-Aramova, Peeva]
  \label{thm:pbs}
  Let \(S=K[x_1,\dots,x_r]\), and suppose that \(I\) is a stable
  monomial ideal in \(S\). 
    Put 
    \begin{displaymath}
      \begin{split}
        b_{i,d} &= \# \setsuchas{x \in G(I)}{\max(x)=i, \, \tdeg{x}=d} \\
        b_i &= \# \setsuchas{x \in G(I)}{\max(x)=i} 
      \end{split}
    \end{displaymath}
    Then, for \(R=S/I\), we have that 
    \begin{align}
      P(\Tor_{*}^{R}(K,K),t) &= 
      \frac{(1+t)^r}{1- t^2 \sum_{i=1}^r (1+t)^{(i-1)} \sum_j b_i } \\
      P(\Tor_{*}^{R}(K,K),t,u) &= 
      \frac{(1+t)^r}{1- t^2 \sum_{i=1}^r (1+tu)^{(i-1)} \sum_j b_{i,j} u^{j}}
  \end{align}
\end{theorem}

Specialising to the case of \(A_n\), we obtain:
\begin{corr}
  \label{corr:pb}
  Let \(n \in \Nat^+\), and let \(r=r(n)\). Regard \(A_n\) as a
  naturally graded \(K\)-algebra, with 
  each \(x_i\) given weight 1, and regard \(K\) as a cyclic \(A\)-module. Then 
  \begin{align}
    P(\Tor_{*}^{A_n}(K,K),t) &= 
    \frac{(1+t)^r}{1- t^2 \sum_{i=1}^r (1+t)^{(i-1)} C_{n,r-i+1} }  \\
    P(\Tor_{*}^{A_n}(K,K),t,u) &= 
    \frac{(1+ut)^r}{1- t^2 \left(\sum_{i=1}^r \left((1+tu)^{(i-1)}
    \sum_j C_{n,r-i+1,j} u^{j}\right)\right)} 
  \end{align}
\end{corr}
Part~\ref{enum:PB}   of the main theorem is now proved.

\begin{example}
  We consider the case \(n=5\), then \(r=r(n)=3\), so
  \(S=K[x_1,x_2,x_3]\) and \(I=I_5 = (x_1^3,\, x_1x_2,\,
  x_1x_3,\, x_2^2, x_2x_3, x_3^2)\). We get that \(C_{5,1}=3,\, C_{5,2}=2,
  C_{5,3}=1\). According to our formulas\footnote{Here, we have used
    the abbreviation  \(P_I^S(t) = P(\Tor_{*}^{S}(I,K),t)\), we will
    also write \(P_K^{S/I}(t) = P(\Tor_{*}^{S/I}(K,K),t)\) et cetera.}
  we  have
  \begin{align*}
  P_I^S(t) &= 1 + 2(1+t) + 3(1+t)^2 =  6  + 8t + 3t^2 \\
  P_K^{S/I} &= \frac{(1+t)^r}{1-t^2 P_I^S(t)} = \frac{1}{1-3t} 
  \end{align*}
  When we consider the grading by total degree, 
  we have that \(C_{5,1,2}=2, \, C_{5,1,3}=1, \, C_{5,2,2} = 2, \,
  C_{5,3,2}=1\). Hence, our formulas yield
  \begin{align*}
    P_I^S(t,u) &= u^2 + 2u^2(1+t) + (2u^2+u^3)(1+t)^2 \\&=
    5{u}^{2}+{u}^{3} + (6{u}^{2}+2{u}^{3})t + \left
      (2\,{u}^{2}+{u}^{3}\right ){t}^{2} \\
    P_K^{S/I}(t,u) &= -{\frac {1+tu}{{u}^{3}{t}^{2}+2\,{t}^{2}{u}^{2}+2\,tu-1}}
  \end{align*}
\end{example}

We list the first few Poincar{\'e}-Betti series
\(P(\Tor_{*}^{A_n}(K,K),t,u)\) in table \ref{fig:betti}.

\begin{conjecture}
\(P(\Tor_{*}^{A_n}(K,K),t) = -\frac{(1+t)^{\ell_1(n)}}{q_n(t)}\), 
\(q_n(t) = \sum_{i=0}^{\ell_2(n)} h_i(n) t^i\), with 
\begin{enumerate}
\item \(q_n(-1) \neq 0\),
\item \(\ell_1(n)\) is the number of odd primes \(p\) such that \(p^2 \le n\), 
\item \(\ell_2(n) = \ell_1(n)+1\),
\item \(h_0(n) = -1\),
\item \(h_1(n) = r(n) - \ell_1(n)\),
\item \(h_{\ell_2(n)}(n) = C_{n,1} = \lceil n/2 \rceil\).
\end{enumerate}
\end{conjecture}

\begin{figure}[p]
  \begin{center}
    \leavevmode
\begin{displaymath}
\begin{array}{|r|c|c|}
  \hline
  n & Graded & Non-graded \\ \hline
2 &-\left (tu-1\right )^{-1}
& -\left (t-1\right )^{-1}
\\
3 &
-\left (2\,tu-1\right )^{-1}
&
-\left (2\,t-1\right )^{-1}
\\
4 &
-{\frac {1+tu}{\left ({u}^{3}+{u}^{2}\right ){t}^{2}+tu-1}}
&
-\left (2\,t-1\right )^{-1}
\\
5 &
-{\frac {1+tu}{\left ({u}^{3}+2\,{u}^{2}\right ){t}^{2}+2\,tu-1}}
&
-\left (3\,t-1\right )^{-1}
\\
6 &
-{\frac {1+tu}{\left (2\,{u}^{3}+{u}^{2}\right ){t}^{2}+2\,tu-1}}
&
-\left (3\,t-1\right )^{-1}
\\
7 &
-{\frac {1+tu}{\left (2\,{u}^{3}+2\,{u}^{2}\right ){t}^{2}+3\,tu-1}}
&
-\left (4\,t-1\right )^{-1}
\\
8 &
-{\frac {1+tu}{\left ({u}^{4}+{u}^{3}+2\,{u}^{2}\right ){t}^{2}+3\,tu-
1}}
&
-\left (4\,t-1\right )^{-1}
\\
9 &
-{\frac {1+2\,tu+{t}^{2}{u}^{2}}{\left ({u}^{5}+2\,{u}^{4}+2\,{u}^{3}
\right ){t}^{3}+\left ({u}^{4}+3\,{u}^{3}+4\,{u}^{2}\right ){t}^{2}+2
\,tu-1}}
&
-{\frac {1+t}{5\,{t}^{2}+3\,t-1}}
\\
10 &
-{\frac {1+2\,tu+{t}^{2}{u}^{2}}{\left ({u}^{5}+3\,{u}^{4}+{u}^{3}
\right ){t}^{3}+\left ({u}^{4}+4\,{u}^{3}+3\,{u}^{2}\right ){t}^{2}+2
\,tu-1}}
&
-{\frac {1+t}{5\,{t}^{2}+3\,t-1}}
\\
11 &
-{\frac {1+2\,tu+{t}^{2}{u}^{2}}{\left ({u}^{5}+3\,{u}^{4}+2\,{u}^{3}
\right ){t}^{3}+\left ({u}^{4}+4\,{u}^{3}+5\,{u}^{2}\right ){t}^{2}+3
\,tu-1}}
&
-{\frac {1+t}{6\,{t}^{2}+4\,t-1}}
\\
12 &
-{\frac {1+2\,tu+{t}^{2}{u}^{2}}{\left (2\,{u}^{5}+2\,{u}^{4}+2\,{u}^{
3}\right ){t}^{3}+\left (2\,{u}^{4}+3\,{u}^{3}+5\,{u}^{2}\right ){t}^{
2}+3\,tu-1}}
&
-{\frac {1+t}{6\,{t}^{2}+4\,t-1}}
\\
13 &
-{\frac {1+2\,tu+{t}^{2}{u}^{2}}{\left (2\,{u}^{5}+2\,{u}^{4}+3\,{u}^{
3}\right ){t}^{3}+\left (2\,{u}^{4}+3\,{u}^{3}+7\,{u}^{2}\right ){t}^{
2}+4\,tu-1}}
&
-{\frac {1+t}{7\,{t}^{2}+5\,t-1}}
\\
14 &
-{\frac {1+2\,tu+{t}^{2}{u}^{2}}{\left (2\,{u}^{5}+3\,{u}^{4}+2\,{u}^{
3}\right ){t}^{3}+\left (2\,{u}^{4}+4\,{u}^{3}+6\,{u}^{2}\right ){t}^{
2}+4\,tu-1}}
&
-{\frac {1+t}{7\,{t}^{2}+5\,t-1}}
\\
15 &
-{\frac {1+2\,tu+{t}^{2}{u}^{2}}{\left (2\,{u}^{5}+4\,{u}^{4}+2\,{u}^{
3}\right ){t}^{3}+\left (2\,{u}^{4}+6\,{u}^{3}+5\,{u}^{2}\right ){t}^{
2}+4\,tu-1}}
&
-{\frac {1+t}{8\,{t}^{2}+5\,t-1}}
\\
16 &
-{\frac {1+2\,tu+{t}^{2}{u}^{2}}{\left ({u}^{6}+{u}^{5}+4\,{u}^{4}+2\,
{u}^{3}\right ){t}^{3}+\left ({u}^{5}+{u}^{4}+6\,{u}^{3}+5\,{u}^{2}
\right ){t}^{2}+4\,tu-1}}
&
-{\frac {1+t}{8\,{t}^{2}+5\,t-1}}
\\
17 &
-{\frac {1+2\,tu+{t}^{2}{u}^{2}}{\left ({u}^{6}+{u}^{5}+4\,{u}^{4}+3\,
{u}^{3}\right ){t}^{3}+\left ({u}^{5}+{u}^{4}+6\,{u}^{3}+7\,{u}^{2}
\right ){t}^{2}+5\,tu-1}}
&
-{\frac {1+t}{9\,{t}^{2}+6\,t-1}}
\\
18 &
-{\frac {1+2\,tu+{t}^{2}{u}^{2}}{\left ({u}^{6}+2\,{u}^{5}+3\,{u}^{4}+
3\,{u}^{3}\right ){t}^{3}+\left ({u}^{5}+2\,{u}^{4}+5\,{u}^{3}+7\,{u}^
{2}\right ){t}^{2}+5\,tu-1}}
&
-{\frac {1+t}{9\,{t}^{2}+6\,t-1}}
\\
19 &
-{\frac {1+2\,tu+{t}^{2}{u}^{2}}{\left ({u}^{6}+2\,{u}^{5}+3\,{u}^{4}+
4\,{u}^{3}\right ){t}^{3}+\left ({u}^{5}+2\,{u}^{4}+5\,{u}^{3}+9\,{u}^
{2}\right ){t}^{2}+6\,tu-1}}
&
-{\frac {1+t}{10\,{t}^{2}+7\,t-1}}
\\
20 &
-{\frac {1+2\,tu+{t}^{2}{u}^{2}}{\left ({u}^{6}+3\,{u}^{5}+2\,{u}^{4}+
4\,{u}^{3}\right ){t}^{3}+\left ({u}^{5}+3\,{u}^{4}+4\,{u}^{3}+9\,{u}^
{2}\right ){t}^{2}+6\,tu-1}}
&
-{\frac {1+t}{10\,{t}^{2}+7\,t-1}}
\\
21 &
-{\frac {\left (1+tu\right )^{2}}{{t}^{3}{u}^{6}+3\,{u}^{5}{t}^{3}+{t}
^{2}{u}^{5}+3\,{t}^{3}{u}^{4}+3\,{u}^{4}{t}^{2}+4\,{t}^{3}{u}^{3}+6\,{
t}^{2}{u}^{3}+8\,{t}^{2}{u}^{2}+6\,tu-1}}
&
-{\frac {1+t}{11\,{t}^{2}+7\,t-1}}
\\
22 &
-{\frac {\left (1+tu\right )^{2}}{{t}^{3}{u}^{6}+3\,{u}^{5}{t}^{3}+{t}
^{2}{u}^{5}+4\,{t}^{3}{u}^{4}+3\,{u}^{4}{t}^{2}+3\,{t}^{3}{u}^{3}+7\,{
t}^{2}{u}^{3}+7\,{t}^{2}{u}^{2}+6\,tu-1}}
&
-{\frac {1+t}{11\,{t}^{2}+7\,t-1}}
\\
23 &
-{\frac {\left (1+tu\right )^{2}}{{t}^{3}{u}^{6}+3\,{u}^{5}{t}^{3}+{t}
^{2}{u}^{5}+4\,{t}^{3}{u}^{4}+3\,{u}^{4}{t}^{2}+4\,{t}^{3}{u}^{3}+7\,{
t}^{2}{u}^{3}+9\,{t}^{2}{u}^{2}+7\,tu-1}}
&
-{\frac {1+t}{12\,{t}^{2}+8\,t-1}}
\\
24 &
-{\frac {\left (1+tu\right )^{2}}{2\,{t}^{3}{u}^{6}+2\,{u}^{5}{t}^{3}+
2\,{t}^{2}{u}^{5}+4\,{t}^{3}{u}^{4}+2\,{u}^{4}{t}^{2}+4\,{t}^{3}{u}^{3
}+7\,{t}^{2}{u}^{3}+9\,{t}^{2}{u}^{2}+7\,tu-1}}
&
-{\frac {1+t}{12\,{t}^{2}+8\,t-1}}
\\
25 &
-{\frac {\left (1+tu\right )^{3}}{q(t,u)}}
&
-{\frac {\left (1+t\right )^{2}}{13\,{t}^{3}+22\,{t}^{2}+7\,t-1}}
\\
\hline
\end{array}  
\end{displaymath}

    \caption{Graded and non-graded Poincar\'e-Betti series of the
  minimal free resolution of \(K\) over \(A_n\).} 
    \label{fig:betti}
  \end{center}
\end{figure}

\end{section}

  \begin{section}{Acknowledgements}
    I am indebted to Johan Andersson for suggesting the idea of studying the
    homological properties of the truncations 
    \(\Gamma_n\). I thank the referee for suggesting a simplified proof of
    parts of 
    Theorem~\ref{thm:bc}. 
  \end{section}

\bibliographystyle{plain}

\end{document}